\documentclass[12pt]{amsart}
\usepackage{amssymb}
\numberwithin{equation}{section}
\theoremstyle{plain}
\newtheorem{theorem}{Theorem}[section]

\newtheorem{cor}[theorem]{Corollary}
\newtheorem{lem}{Lemma}[section]

\theoremstyle{definition}

\theoremstyle{remark}
\newtheorem{rem}{Remark}[section]

\def\R{{\mathbb R}}
\def\N{{\mathbb N}}
\def\Z{{\mathbb Z}}
\def\p#1{{\left({#1}\right)}}
\def\b#1{{\left\{{#1}\right\}}}
\def\br#1{{\left[{#1}\right]}}
\def\jp#1{{\left\langle{#1}\right\rangle}}
\def\n#1{{\left\|{#1}\right\|}}
\def\abs#1{{\left|{#1}\right|}}

\def\supp{\operatorname{supp}}
\def\Rnx{{\mathbb R}^n_x}
\def\Rnxi{{\mathbb R}^n_\xi}
\def\Rone{{\mathbb R}}
\def\Rn{{{\mathbb R}^n}}
\def\Scal{\mathcal{S}}
\def\Fcal{\mathcal{F}}
\def\Nat{{\mathbb N}}

\title[]{Global $L^2$-boundedness theorems for a class
   of Fourier integral operators}

\author[]{Michael Ruzhansky and Mitsuru Sugimoto}
\address{
  Michael Ruzhansky:
  \endgraf
  Department of Mathematics
  \endgraf
  Imperial College of Science, Technology and Medicine
  \endgraf
  180 Queen's Gate, London SW7 2BZ, UK
  \endgraf
  {\it E-mail address} {\rm ruzh@ic.ac.uk}
  \endgraf
  \medskip
  Mitsuru Sugimoto:
  \endgraf
  Department of Mathematics, Graduate School of Science
  \endgraf
  Osaka University
  \endgraf
  Machikaneyama-cho 1-16, Toyonaka, Osaka 560-0043, Japan
  \endgraf
  {\it E-mail address} {\rm sugimoto@math.wani.osaka-u.jp}
  }
\begin{document}
\begin{abstract}
The local $L^2$-mapping property of Fourier integral operators has
been established in H\"ormander \cite{H} and in Eskin \cite{E}. In
this paper, we treat the global $L^2$-boundedness for a class of
operators that appears naturally in many problems. As a consequence, 
we will improve
known global results for several classes of pseudo-differential and
Fourier integral operators, as well as extend previous results of
Asada and Fujiwara \cite{AF} or Kumano-go \cite{Ku}. As an
application, we show a global smoothing estimate to generalized
Schr\"odinger equations which extends the results of Ben-Artzi and
Devinatz \cite{BD}, Walther \cite{Wa}, and \cite{Wa2}.
\end{abstract}

\maketitle
%\newpage
\section{Introduction}
We consider (Fourier integral) operators, which can be globally
written in the form
\begin{equation}\label{eq:eq1}
      Tu(x)=\int_{\R^n}\int_{\R^n}
      e^{i\phi(x,y,\xi)}a(x,y,\xi)u(y)\, dy d\xi
\quad
(x\in\R^n),
\end{equation}
where $a(x,y,\xi)$ is an amplitude function and $\phi(x,y,\xi)$ is
 a real phase function of the form
\[
\phi(x,y,\xi)=x\cdot\xi+\varphi(y,\xi).
\]
Note that, by the equivalence of phase function theorem, Fourier
integral operators with the local graph condition can always be
written in this form locally.
Although, due to the
nontriviality of the Maslov cohomology class, globally defined
Fourier integral operators can not be written in this form with a
globally defined phase $\phi$, it is nevertheless convenient to
still call them Fourier integral operators.
It will be clear below
how such operators naturally arise in global smoothing problems if
we use an adaptation of the Egorov theorem.
\par
Local $L^2$ mapping property of the operator \eqref{eq:eq1} has
been established in H\"ormander \cite{H} and in Eskin \cite{E}.
One of the aims of this paper is to establish the global
$L^2$-boundedness properties of operators (\ref{eq:eq1}). 
Analogous properties can be then easily obtained for adjoint operators
as well.

We will
try to make as few assumptions as possible in the spirit of global
$L^2$-estimates for pseudo-differential operators (see
Calder\'on--Vaillancourt \cite{CV}, Childs \cite{Chil},
Coifman-Meyer \cite{CM}, Cordes \cite{Co}). In fact, our Corollary
\ref{Cor2.4} will not only extend these $L^2$-boundedness results
to more general operators (\ref{eq:eq1}), but will also reduce the
number of assumptions on the amplitude in the case of
pseudo-differential operators, compared to the above mentioned papers
(see also Sugimoto \cite{Su}).
Global $L^2$-boundedness of operators (\ref{eq:eq1}) has been
previously studied by Asada-Fujiwara \cite{AF}, Kumano-Go
\cite{Ku}. However, there
one had to make a quite restrictive and not always natural
assumption on the boundedness of $\partial_\xi\partial_\xi\phi,$
which fails in many important cases.
In Coriasco \cite{Cor} and Boggiato-Buzano-Rodino
\cite{BBR} such results are applied to obtain global estimates of
solutions to some classes of hyperbolic equations.
Here again one requires quite strong decay properties of derivatives
of both phase and amplitude. We will remove all these assumptions
and will give general $L^2$-estimates.
In fact, in global estimates of Section 2 we will actually impose only 
a finite number of conditions on the phase and the amplitude, compared to
infinitely many in the above mentioned papers. 
\par
As a consequence of our $L^2$-estimates, we can treat canonical
transforms.
Operators that appear
there are of the form (\ref{eq:eq1}) with phase function
\begin{equation}\label{eq:eq2}
\phi(x,y,\xi)=x\cdot\xi-y\cdot p(\xi)\frac{\nabla p(\xi)}{\abs{\nabla p(\xi)}},
\end{equation}
where
$p(\xi)$ is a positively homogeneous function of degree $1$.
If we take $p(\xi)=|\xi|$, then we have $\phi(x,y,\xi)=x\cdot\xi-y\cdot\xi$,
and the operator $T$ defined by \eqref{eq:eq1}
is a pseudo-differential operator.
Furthermore, the operator $T$ with general \eqref{eq:eq2} is used to transform
the Fourier multiplier
\[
L_p=p(D_x)^2=F^{-1}_\xi p(\xi)^2F_x
\]
to the
Laplacian $-\triangle$, where $F_x$ ($F_\xi^{-1}$ {\it resp.})
denotes the (inverse {\it resp.})
Fourier transform.
In fact, we have a relation
\[
T\cdot(-\triangle)\cdot T^{-1}=L_p
\]
under a certain condition on $p(\xi)$ if we take $1$ as the amplitude
function $a(x,y,\xi)$ (see Section 4).
The $L^2$-property of the Laplacian is well known in various
situations. Our objective is to know the $L^2$-property of the
operator $T$, so that we can extract the $L^2$-property of the
operator $L_p$ from that of the Laplacian.
This approach allows to give a general treatment of several
smoothing problems, including those treated by e.g.
Ben-Artzi and Klainerman \cite{BK}, Simon \cite{Si},
Kato and Yajima \cite{KY}, or Walther \cite{Wa}.
\par
We should mention here that the global $L^2$-boundedness
with example (\ref{eq:eq2})
is not covered by previous results,
for example, Asada and Fujiwara \cite{AF},
and Kumano-go \cite{Ku}.
The result of \cite{AF} is motivated by the construction of fundamental
solution of Schr\"odinger equation in the way of Feynman's path integral,
and it requires the boundedness of all the derivatives of entries of
the matrix
\[
\begin{pmatrix}
\partial_x\partial_y\phi & \partial_x\partial_\xi\phi
\\
\partial_\xi\partial_y\phi & \partial_\xi\partial_\xi\phi
\end{pmatrix}.
\]
For the details, see \cite{AF} and references cited there.
With our example (\ref{eq:eq2}), the boundedness of the entries of
$\partial_\xi\partial_\xi\phi$ fails.
On the other hand, the result of \cite{Ku} is used to construct the
fundamental solution of hyperbolic equations, and it requires that
$J(y,\xi)=\phi(x,y,\xi)-(x-y)\cdot\xi$ satisfies
\[
\abs{\partial_y^\alpha\partial_\xi^\beta J(y,\xi)}
\leq C_{\alpha\beta}\p{1+|\xi|}^{1-|\beta|}
\]
for all $\alpha$ and $\beta$.
Our example (\ref{eq:eq2}) does not satisfy these estimates with $\alpha=0$.
\par
In this paper, we develop a new $L^2$-theory which does not
require these decay assumptions. In particular, it includes the
case of example \eqref{eq:eq2}. For $m\in\R$, let $L^2_m(\R^n)$ be
the set of functions $f$ such that the norm
\[
\n{f}_{L^2_m(\R^n)}
=\p{\int_{\R^n}\abs{\langle x\rangle^m f(x)}^2\,dx}^{1/2};
\quad
\langle x\rangle^m=\p{1+|x|^2}^{m/2}
\]
is finite.
The following is a simplified version of our main result
(Theorem \ref{Th3.1}) which is expected to have many applications:
\medskip
\begin{theorem}\label{mainth}
Let the operator $T$ be defined by \eqref{eq:eq1},
where
$\varphi(y,\xi)\in C^\infty\p{\R^n_y\times\R^n_\xi}$ is a real-valued function,
and $a(x,y,\xi)\in C^\infty\p{\R^n_x\times\R^n_y\times\R^n_\xi}$.
Assume that
\[
\abs{\det\partial_y\partial_{\xi}\varphi(y,\xi)}\geq C>0,
\]
and all the derivatives of entries of $\partial_y\partial_{\xi}\varphi$
are bounded.
Also assume that
\begin{align*}
&\abs{\partial_\xi^\alpha\varphi(y,\xi)} \leq
C_{\alpha}\langle y\rangle\quad
\text{for all $|\alpha|\geq1$},
\\
&\abs{\partial_x^\alpha\partial_y^\beta\partial_\xi^\gamma a(x,y,\xi)}
\leq C_{\alpha\beta\gamma}\langle x\rangle^{-|\alpha|}
\quad
\text{for all $\alpha$, $\beta$, and $\gamma$}.
\end{align*}
Then $T$ is bounded on $L^2_{m}(\R^n)$ for any $m\in\R$.
\end{theorem}
\medskip
This theorem says that, if amplitude functions $a(x,y,\xi)$ have
some decaying properties with respect to $x$,
we do not need the boundedness of $\partial_\xi\partial_\xi\phi$
for the $L^2$-boundedness, as required in \cite{AF},
and we can have weighted estimates, as well.
(The same is true when both phase and amplitude functions have
some decaying properties with respect to $y$.
See Theorem \ref{Th3.1}.)
\par
We explain the plan of this paper.
In Section 2, we show the global $L^2$-boundedness of a class
of oscillatory integral operators, which
generalizes a standard local result explained in Stein \cite{St}.
By using it, we prove various type of the $L^2$-boundedness of Fourier
integral operators.
Some of them are extension of previous results on the
$L^2$-boundedness of pseudo-differential operators with non-regular
symbols.
It is worth mentioning that, in general, we do not necessarily need the
standard homogeneity assumption for the phase function in the
frequency variable.
In addition, we impose the boundedness condition on
only a finite number of the derivatives of phase functions,
while infinitely many in \cite{AF} and \cite{Ku}.
\par
In Section 3, we state and prove our main result Theorem
\ref{Th3.1}. We remark that it (together with Theorem \ref{Th2.5})
substantially weaken the assumptions for the $L^2$-boundedness of
SG pseudo-differential (as in Cordes \cite{Co1}) and SG Fourier
integral operators (as in Coriasco \cite{Cor}). These operators
are used to handle the SG hyperbolic partial differential
equations (roughly speaking, certain equations with coefficients of
polynomial growth). The class of symbols $SG^{m_1,m_2}$ is defined
as a space of smooth functions $a=a(y,\xi)\in
C^\infty(\R^n_y\times\R^n_\xi)$ satisfying the estimate
\[
|\partial_y^\beta\partial_\xi^\gamma a(y,\xi)|\leq C_{\beta\gamma}
 \langle y \rangle^{m_1-|\beta|} \langle \xi \rangle^{m_2-|\gamma|}
  \;\;\; \text{for all $\beta$ and $\gamma$}.
\]
SG Fourier integral operators are operators of the form
\[
Tu(x)=\int_{\R^n} e^{i(x\cdot\xi+\varphi(y,\xi))}a(y,\xi)u(y)d\xi dy
\]
(or its adjoint),
where $a\in SG^{m_1,m_2}$ and $\varphi\in SG^{1,1}$, which also satisfies
\[
   C_1\langle y \rangle \leq \langle \partial_\xi\varphi \rangle \leq
   C_2\langle y \rangle, \;\;\;
   C_1\langle \xi \rangle \leq \langle \partial_y\varphi \rangle \leq
   C_2\langle \xi \rangle,
\]
for some $C_1, C_2>0$. A result in \cite{Co1} for SG
pseudo-differential and its extension in \cite{Cor} for SG Fourier
integral operators states that under these assumptions on the
phase $\phi$, and for $a\in SG^{0,0}$, the corresponding operator
$T$ is bounded on $L^2(\R^n)$. Without going much into detail, let
us mention here that statements of our results replace the strong
decay assumptions $\phi\in SG^{1,1}$, $a\in SG^{0,0}$, by (a
finite number of) boundedness conditions, for $T$ to be still
bounded in $L^2(\R^n)$.
\par
In Section 4, we exhibit an example of how to use our main result.
We mainly focus on the problem of global smoothing property
of generalized Schr\"odinger equations
\begin{equation}\label{eq:Q}
\left\{
\begin{array}{l}
\p{i\partial_t+Q(D)}u(t,x)=0, \\
u(0,x)=f(x).
\end{array}
\right.
\end{equation}
Ben-Artzi and Devinatz \cite{BD} showed a global smoothing estimate
to equation \eqref{eq:Q}, where the symbol $Q(\xi)$ of $Q(D)$ is
a real polynomial of principal type.
Walther \cite{Wa2} consider the case of radially symmetric $Q(\xi)$.
By using our result Theorem \ref{Th3.1}, we can treat more general case
(see Theorem \ref{Th4.2}).
More refined applications to this subject
will be shown in our forthcoming paper
\cite{RS3}.
In subsequent work \cite{RS2}, we will establish properties of operators
(\ref{eq:eq1}) in weighted Sobolev spaces, which will have several
further applications of these results to hyperbolic equations as
well as global canonical transforms. 
\par
%\newpage
\bigskip
\section{Global $L^2$-estimates}
First of all,
we confirm a basic result on the
$L^2$-boundedness of a class of oscillatory integral operators,
based on the argument of Fujiwara \cite{F},
which is a global version of a proposition in Stein
\cite[p.377]{St}.
Here and hereafter, the capital $C$ (sometimes
with some suffices) always denotes a positive constant which may
differ on each occasion.
\medskip
\begin{theorem}\label{Th2.1}
Let the operator $I_\varphi$ be defined by
\begin{equation}\label{OsI}
I_\varphi u(x)=\int_{\R^n}e^{i\varphi(x,y)}a(x,y)u(y)\,dy,
\end{equation}
where $a(x,y)\in C^\infty\p{\R^n_x\times\R^n_y}$, and
$\varphi(x,y)\in C^\infty\p{\R^n_x\times\R^n_y}$ is a
real-valued function.
Assume that
\[
\abs{\partial_x^\alpha\partial_y^\beta a(x,y)}
\leq C_{\alpha\beta},
\]
for $|\alpha|, |\beta|\leq 2n+1$.
Also assume that,  on $\supp a(x,y)$,
\[
\abs{\det\partial_x\partial_y\varphi(x,y)}\geq C>0
\]
and each entry $h(x,y)$ of
the matrix $\partial_x\partial_y\varphi(x,y)$ satisfies
\[
\abs{\partial_x^\alpha h(x,y)}\leq C_{\alpha},
\qquad
\abs{\partial_y^\beta h(x,y)}\leq C_{\beta}
\]
for $|\alpha|,|\beta|\leq 2n+1$.
Then the operator $I_\varphi$ is $L^2(\R^n)$-bounded, and satisfies
\[
\n{I_\varphi}_{L^2\to L^2}
\leq C\sup_{|\alpha|,|\beta|\leq 2n+1}
 \n{\partial_x^\alpha\partial_y^\beta a(x,y)}_{L^\infty\p{\R^n_x\times\R^n_y}}.
\]
\end{theorem}
\medskip
\begin{proof}
Let $g\in C^\infty_0(\R^n)$ be a real-valued positive
function such that $\{g_k(x)\}_{k\in\Z^n}$, where $g_k(x)=g(x-k)$,
forms a partition of unity.
We decompose the operator $I_\varphi$ as
\[
I_\varphi=\sum_{(j,k)\in\Z^n\times\Z^n}I_{(j,k)},
\]
where $I_{(j,k)}=g_jI_\varphi g_k$, that is,
\[
I_{(j,k)}u(x)=g_j(x)\int e^{i\varphi(x,z)}a(x,z)g_k(z)u(z)\,dz.
\]
We denote the adjoint of $I_{(j,k)}$ by $I_{(j,k)}^*$,
that is,
\[
I_{(j,k)}^*u(z)=g_k(z)\int e^{-i\varphi(y,z)}
\overline{a(y,z)}g_j(y)u(y)\,dy.
\]
Then we have
\[
I_{(j,k)}I_{(l,m)}^*u(x)=\int K_{(j,k),(l,m)}(x,y)u(y)\,dy,
\]
where
\[
K_{(j,k),(l,m)}(x,y)=g_j(x)g_l(y)
  \int e^{i\p{\varphi(x,z)-\varphi(y,z)}}
   a(x,z)\overline{a(y,z)}
  g_k(z)g_m(z)\,dz.
\]
By integration by parts, we have
\begin{align*}
&\int e^{i\p{\varphi(x,z)-\varphi(y,z)}}   a(x,z)\overline{a(y,z)}
g_k(z)g_m(z)\,dz
\\
=&\int e^{i\p{\varphi(x,z)-\varphi(y,z)}}
L^{2n+1}\p{   a(x,z)\overline{a(y,z)}
g_k(z)g_m(z)}\,dz,
\end{align*}
where $L$ is the transpose of the operator
\[
^tL=\frac 1i
\frac
{\partial_z\varphi(x,z)-\partial_z\varphi(y,z)}
{\abs{\partial_z\varphi(x,z)-\partial_z\varphi(y,z)}^2}
\cdot\partial_z
\]
From the assumptions, and using that
\[
\partial_z\phi(x,z)-\partial_z\phi(y,z)=\partial_x\partial_z\phi(w,z)(x-y)
\]
for some $w$, we obtain
\[
\abs{\partial_z\varphi(x,z)-\partial_z\varphi(y,z)}\geq C|x-y|
\]
and
\[
\abs{\partial_z^\beta\varphi(x,z)-\partial_z^\beta\varphi(y,z)}
\leq C_\beta|x-y|
\]
for $1\leq|\beta|\leq 2n+2$.
Hence, we have
\[
|K_{(j,k),(l,m)}(x,y)|
\leq CA^2\frac {g_j(x)g_l(y)}{1+|x-y|^{2n+1}}h(k-m),
\]
where $h\in C_0^\infty(\R^n)$ is a positive function
($h(x)=\int g(z-x)g(z)\,dz$),
and
\[
A=\sup_{|\alpha|,|\beta|\leq 2n+1}
 \n{\partial_x^\alpha\partial_y^\beta a}_{L^\infty\p{\R^n_x\times\R^n_y}}.
\]
Then we have
\begin{align*}
&\sup_{x}\int |K_{(j,k),(l,m)}(x,y)|\,dy
    \leq CA^2 \frac {h(k-m)}{1+|j-l|^{2n+1}},
\\
&\sup_{y}\int |K_{(j,k),(l,m)}(x,y)|\,dx
    \leq CA^2 \frac {h(k-m)}{1+|j-l|^{2n+1}},
\end{align*}
which implies
\[
\n{I_{(j,k)}I_{(l,m)}^*}_{L^2\to L^2}
    \leq CA^2 \frac {h(k-m)}{1+|j-l|^{2n+1}}.
\]
Here we have used the following lemma
(see Stein \cite[p.284]{St}):
\medskip
\begin{lem}\label{Lem2.1}
Suppose $S$ is given by
\[
(Sf)(x)=\int s(x,y)f(y)\,dy,
\]
where the kernel $s(x,y)$ satisfies
\[
\sup_x\int|s(x,y)|\,dy\leq 1,\qquad
\sup_y\int|s(x,y)|\,dx\leq 1.
\]
Then $\n{S}_{L^2\to L^2}\leq 1$.
\end{lem}
\medskip
By the same discussion, we have
\[
\n{I_{(j,k)}^*I_{(l,m)}}_{L^2\to L^2}
    \leq CA^2 \frac {h(j-l)}{1+|k-m|^{2n+1}}.
\]
Then we have
\[
\n{I_{(j,k)}I_{(l,m)}^*}_{L^2\to L^2},\,
\n{I_{(j,k)}^*I_{(l,m)}}_{L^2\to L^2}
    \leq CA^2 \b{\gamma\p{j-l,k-m}}^2,
\]
where
\[
\gamma\p{j_1,j_2}
=\sqrt{\b{\frac {h(j_2)}{1+|j_1|^{2n+1}}+\frac {h(j_1)}{1+|j_2|^{2n+1}}}}
\]
and it satisfies the estimate
\[
\sum_{(j_1,j_2)\in\Z^n\times\Z^n}
\gamma\p{j_1,j_2}<\infty.
\]
We have the desired result, by the following Cotlar's lemma
(see Calder\'on and Vaillancourt \cite{CV},
Stein \cite[Chapter VII, Section 2]{St}):
\medskip
\begin{lem}\label{Lem2.2}
Assume a family of $L^2$-bounded operators $\{T_j\}_{j\in\Z^r}$
and positive constants $\{\gamma(j)\}_{j\in\Z^r}$ satisfy
\[
\n{T_i^*T_j}_{L^2\to L^2}\leq \b{\gamma(i-j)}^2,\qquad
\n{T_iT_j^*}_{L^2\to L^2}\leq \b{\gamma(i-j)}^2,\qquad
\]
and
\[
M=\sum_{j\in\Z^r}\gamma(j)<\infty.
\]
Then the operator
\[
T=\sum_{j\in\Z^r}T_j
\]
satisfies
\[
\n{T}_{L^2\to L^2}\leq M.
\]
\end{lem}
\medskip
\end{proof}
\medskip
By using Theorem \ref{Th2.1} on oscillatory integral operators \eqref{OsI},
we can easily show the $L^2$-boundedness of
Fourier integral operators of special forms.
Let us begin with the case when
the amplitude $a(x,y,\xi)$ is independent of the variable $y$.
\medskip
\begin{theorem}\label{Th2.2}
Let the operator $T$ be defined by
\begin{equation}\label{T1}
 Tu(x)
 =\int_{\R^n}\int_{\R^n} e^{i(x\cdot\xi+\varphi(y,\xi))}a(x,\xi)u(y) dy d\xi,
\end{equation}
where $a(x,\xi)\in C^\infty\p{\R^n_x\times\R^n_\xi}$ and
$\varphi(y,\xi)\in C^\infty\p{\R^n_y\times\R^n_\xi}$.
Assume that the pseudo-differential operators $a(X,D)$ defined by
\[
a(X,D)u(x)
 =(2\pi)^{-n}\int_{\R^n} \int_{\R^n}e^{i(x-y)\cdot\xi}a(x,\xi)u(y) \,dyd\xi
\]
and the oscillatory integral operator $I_\varphi$
defined by
\[
I_\varphi u(\xi)=\int_{\R^n} e^{i\varphi(y,\xi)}\chi_E(\xi)u(y)\,dy
\]
are both $L^2(\R^n)$-bounded,
where $\chi_E$ is the characteristic function of the set
\begin{equation}\label{support}
E=\bigcup_{x\in\R^n}E_x\,;\quad
E_x=\supp a(x,\cdot)\subset\R^n.
\end{equation}
Then $T$ is $L^2(\R^n)$-bounded, and satisfies
\[
\n{T}_{L^2\to L^2}
\leq \p{2\pi}^{n/2}\n{a(X,D)}_{L^2\to L^2}\cdot\n{I_\varphi}_{L^2\to L^2}.
\]
\end{theorem}
\medskip
\begin{proof}
We remark that
$T=\p{2\pi}^na(X,D) F^{-1}I_\varphi$, where $F^{-1}$ is the inverse Fourier
transform.
The $L^2\p{\R^n}$-boundedness of $T$ is obtained from the assumptions
and Plancherel's theorem.
\end{proof}
\medskip
As a corollary, we have the result announced
in Ruzhansky and Sugimoto \cite{RS}.
Now we recall the definition of the Besov space
  $B_{p,q}^{({\bf s,s^\prime})}$ for $0<p,q\leq\infty$ and multi-indices
  $({\bf s,s^\prime})$, where ${\bf s}=(s_1,\ldots,s_N)$
  and ${\bf s^\prime}=(s_1^\prime,\cdots,s_{N^\prime}^\prime)$.
  Let ${\bf n}=(n_1,\ldots,n_N)$,
 ${\bf n^\prime}=(n^\prime_1,\cdots,n^\prime_{N^\prime})$
  be splitting of $\Rnx$ and $\Rnxi$, respectively:
\[
n=n_1+\ldots +n_N=n_1^\prime+\ldots +n^\prime_{N^\prime}.
\]
Then
  $f\in B_{p,q}^{({\bf s, s^\prime})}=
    B_{p,q}^{({\bf s,s^\prime})}({\Rone}^{({\bf n,\bf n^\prime})})$ if
  $f=f(x,\xi)\in\Scal^\prime(\Rone^{2n})$ and
\[
\n{f}_{B_{p,q}^{({\bf s,s^\prime})}}=
   \left\{ \sum_{{\bf j,k}\geq 0}\left( \int_\Rn\int_\Rn
    |2^{{\bf j\cdot s}+{\bf k\cdot s^\prime}}
    \Fcal^{-1}\Phi_{\bf j,k}\Fcal f(x,\xi)|^p
      dx d\xi\right)^{q/p}\right\}^{1/q}<\infty.
\]
Here ${\bf j}=(j_1,\ldots,j_N)$, ${\bf k}=(k_1,\ldots,k_{N^\prime})$,
 $\Fcal$ is the Fourier transform with respect to $(x,\xi)$,
 $\Fcal^{-1}$ is the inverse Fourier transform with respect to the
dual variable $(y,\eta)$,
and
$\Phi_{\bf j,k}=\Phi_{\bf j,k}(y,\eta)=\Theta_{j_1}(y_1)\cdots\Theta_{j_N}(y_N)
   \Theta_{k_1}(\eta_1)\cdots\Theta_{k_{N^\prime}}(\eta_{N^\prime})$.
    Here we split variables $y,\eta\in\Rone^n$ following
    the splitting ${\bf n}$, ${\bf n'}$.
Functions
 $\Theta_{i}(z)\in\Scal$ form the dyadic system
  of the corresponding dimension:
   $\supp\Theta_0\subset\{z;|z|\leq 2\}$,
   $\supp\Theta_i\subset\{z;2^{i-1}\leq |z|\leq 2^{i+1}\}$ for $i\in{\Nat}$,
   $\sum_{i=0}^\infty \Theta_i(z)=1$, and
   $2^{i|\alpha|}|\partial^\alpha\Theta_i(z)|\leq C_\alpha$
for all $i\geq 0$ and all
   $z$.
A natural modification is needed for $p,q=\infty$,
 see \cite{Tr}.
\medskip
 \begin{cor}\label{Cor2.3}
Let $2\leq p\leq\infty$. Let the operator $T$ be defined by
{\rm{(\ref{T1})}}, where $\varphi(y,\xi)\in
C^\infty\p{\R^n_y\times\R^n_\xi}$ is a real-valued function, and
$a(x,\xi)\in B_{p,1}^{(1/2-1/p)({\bf n},{\bf n^\prime})}$.
Assume that, on $\R^n\times E$,
\[
\abs{\det\partial_y\partial_{\xi}\varphi(y,\xi)}\geq C>0
\]
and each entry $h(y,\xi)$ of
the matrix $\partial_y\partial_\xi\varphi(y,\xi)$ satisfies
\[
\abs{\partial_y^\alpha h(y,\xi)}\leq C_{\alpha},
\qquad
\abs{\partial_\xi^\beta h(y,\xi)}\leq C_{\beta}
\]
for $|\alpha|,|\beta|\leq 2n+1$, where $E$
is the set defined by \eqref{support}.
Then $T$ is $L^2(\R^n)$-bounded, and satisfies
\[
\n{Tu}_{L^2(\Rn)}\leq C\n{a(x,\xi)}
      _{B_{p,1}^{(1/2-1/p)({\bf n},{\bf n^\prime})}}
  \n{u}_{L^2(\Rn)}.
\]
\end{cor}
\medskip
\begin{proof}
The $L^2$-boundedness of $T$ follows from Theorems \ref{Th2.1},
\ref{Th2.2}, and the fact
that pseudo-differential operators $a(X,D)$ with
$a(x,\xi)\in B_{p,1}^{(1/2-1/p)({\bf n},{\bf n^\prime})}$
are $L^2$-bounded.
See Sugimoto \cite{Su}.
\end{proof}
\medskip
Corollary \ref{Cor2.3} is rather general but its conditions
 may be hard to check.
On the other hand, conditions of the corollary below can
be checked in various situations.
\medskip
\begin{cor}\label{Cor2.4}
 Let the operator $T$ be defined by {\rm{(\ref{T1})}},
where $\varphi(y,\xi)\in C^\infty\p{\R^n_y\times\R^n_\xi}$
is a real-valued function.
Assume that, on $\R^n\times E$,
\[
\abs{\det\partial_y\partial_{\xi}\varphi(y,\xi)}\geq C>0
\]
and each entry $h(y,\xi)$ of
the matrix $\partial_y\partial_\xi\varphi(y,\xi)$ satisfies
\[
\abs{\partial_y^\alpha h(y,\xi)}\leq C_{\alpha},
\qquad
\abs{\partial_\xi^\beta h(y,\xi)}\leq C_{\beta}
\]
for $|\alpha|,|\beta|\leq 2n+1$, where $E$
is the set defined by \eqref{support}.
Also assume that $a(x,\xi)$ belongs to the symbol class $S_{0,0}^0$
(that is,
 $\partial_x^\alpha\partial_\xi^\beta a(x,\xi)\in L^\infty(\Rnx\times\Rnxi)$
 for all $\alpha$ and $\beta$).
Otherwise assume one of the following conditions:
 \begin{itemize}
  \item[{\rm (1)}] $\partial_x^\alpha\partial_\xi^\beta a(x,\xi)\in
    L^\infty(\Rnx\times\Rnxi)$ for $\alpha,\beta \in\{0,1\}^n$.
  \item[{\rm (2)}] $\partial_x^\alpha\partial_\xi^\beta a(x,\xi)\in
    L^\infty(\Rnx\times\Rnxi)$
    for $|\alpha|,|\beta| \leq [n/2]+1$.
 \item[{\rm (3)}] $\partial_x^\alpha\partial_\xi^\beta a(x,\xi)\in
    L^\infty(\Rnx\times\Rnxi)$
    for $|\alpha|\leq [n/2]+1$, $\beta\in\{0,1\}^n$.
 \item[{\rm (4)}] $\partial_x^\alpha\partial_\xi^\beta a(x,\xi)\in
    L^\infty(\Rnx\times\Rnxi)$
    for $\alpha\in\{0,1\}^n$, $|\beta|\leq [n/2]+1$.
 \item[{\rm (5)}] There exist real numbers $\lambda,\lambda^\prime>n/2$
   such that $(1-\Delta_x)^{\lambda/2} (1-\Delta_\xi)^{\lambda^\prime/2}
   a(x,\xi)\in L^\infty(\Rnx\times\Rnxi)$.
 \item[{\rm (6)}] There exist a real number $\lambda>1/2$ and a constant $C$
  such that $$||\delta_x^\alpha(h)\delta_\xi^\beta(h^\prime) a(x,\xi)||_{
  L^\infty(\Rnx\times\Rnxi)}\leq C\prod_{i,j=1}^n |h_i|^{\alpha_i\lambda}
  |h_j^\prime|^{\beta_j\lambda}$$ holds for all
  $\alpha,\beta\in\{0,1\}^n$ and all
  $h=(h_1,\ldots,h_n), h^\prime=(h_1^\prime,\ldots,h_n^\prime)\in \Rn$.
  Here $\delta_x^\alpha(h)=\delta_{x_1}^{\alpha_1}(h_1)
     \cdots\delta_{x_n}^{\alpha_1}(h_n)$
  is the difference operator, with
\[
  \delta_{x_i}^0(h_i)a(x,\xi)=a(x,\xi),\qquad
  \delta_{x_i}^1(h_i)a(x,\xi)=a(x+h_i e_i,\xi)-a(x,\xi),
\]
 where $e_i$ is the $i$-th
  standard basis vector in $\Rn$. The definition of $\delta_\xi^\beta$ is
 similar.
 \item[{\rm (7)}] There exists a real number $2\leq p<\infty$ such that
    $\partial_x^\alpha\partial_\xi^\beta a(x,\xi)\in
    L^p(\Rnx\times\Rnxi)$
    for $|\alpha|,|\beta|\leq [n(1/2-1/p)]+1$.
 \end{itemize}
  Then $T$ is $L^2(\Rn)$-bounded.
  \end{cor}
\medskip
Corollary \ref{Cor2.4} with $\varphi(y,\xi)=-y\cdot\xi$ is a refined version
of known results on the $L^2$-boundedness of pseudo-differential operators with
non-regular symbols;
(1) with $\alpha,\beta\in\{0,1,2,3\}^n$ is due to Calder\'on and Vaillancourt
 \cite{CV}, (2) and (5) are due to
Cordes \cite{Co}, the difference condition (6) is due to Childs \cite{Chil},
and conditions (3) with $|\alpha|\leq[n/2]+1$, $\beta\in\{0,1,2\}^n$,
(7) with $\alpha\leq[n(1/2-1/p)]+1$, $|\beta|\leq 2n$
are due to Coifman and Meyer \cite{CM}.
\medskip
\begin{proof}
The $L^2$-boundedness under all conditions follows from Corollary \ref{Cor2.3}
 with different choices of $p$ and splitting of the spaces.
In fact, conditions (1)--(6) are obtained with
 $p=\infty$ ($N=N^\prime=1$ in conditions (2), (5);
 $N=N^\prime=n$ in (1), (6); $N=1, N^\prime=n$ in (3), and $N=n, N'=1$ in (4)).
 Condition (7) is obtained from Corollary \ref{Cor2.3} by taking
 the same $p$ and $N=N^\prime=1$.
 For more relations between symbol classes and
 Besov spaces, we refer to Sugimoto \cite{Su} and Triebel \cite{Tr}.
\end{proof}
\medskip
\par
We have a theorem for amplitudes which are independent of the variable $x$,
as well.
\medskip
\begin{theorem}\label{Th2.5}
Let the operator $T$ be defined by
\begin{equation}\label{T2}
 Tu(x)
 =\int_{\R^n}\int_{\R^n} e^{i(x\cdot\xi+\varphi(y,\xi))}a(y,\xi)u(y) dy d\xi,
\end{equation}
where $a(y,\xi)\in C^\infty\p{\R^n_y\times\R^n_\xi}$, and
$\varphi(y,\xi)\in C^\infty\p{\R^n_y\times\R^n_\xi}$ is
a real-valued function.
Assume that
\[
\abs{\partial_y^\alpha\partial_\xi^\beta a(y,\xi)}
\leq C_{\alpha\beta},
\]
for $|\alpha|, |\beta|\leq 2n+1$.
Also assume that, on $\supp a(y,\xi)$,
\[
\abs{\det\partial_y\partial_\xi\varphi(y,\xi)}\geq C>0
\]
and each entry $h(y,\xi)$ of
the matrix $\partial_y\partial_\xi\varphi(y,\xi)$ satisfies
\[
\abs{\partial_y^\alpha h(y,\xi)}\leq C_{\alpha},
\qquad
\abs{\partial_\xi^\beta h(y,\xi)}\leq C_{\beta}
\]
for $|\alpha|,|\beta|\leq 2n+1$.
Then the operator $T$ is $L^2(\R^n)$-bounded, and satisfies
\[
\n{T}_{L^2\to L^2}\leq C
\sup_{|\alpha|,|\beta|\leq 2n+1}
 \n{\partial_y^\alpha\partial_\xi^\beta a(y,\xi)}
  _{L^\infty\p{\R^n_y\times\R^n_\xi}}.
\]
\end{theorem}
\medskip
\begin{proof}
We remark that $T=(2\pi)^nF^{-1}I_\varphi$, where $F^{-1}$ is the inverse
Fourier transform and $I_\varphi$ is the oscillatory integral operator
defined by
\[
I_\varphi u(\xi)=\int e^{i\varphi(y,\xi)}a(y,\xi)u(y)\,dy.
\]
The result is obtained from Theorem \ref{Th2.1} and Plancherel's theorem.
\end{proof}
\medskip
\par
As a corollary of Theorems \ref{Th2.2} and \ref{Th2.5}, we have a result for
amplitudes which are of the product type.
\medskip
\begin{cor}\label{Cor2.6}
Let the operator $T$ be defined by
\begin{align*}
&Tu(x)
 =\int_{\R^n}\int_{\R^n} e^{i(x\cdot\xi+\varphi(y,\xi))}a(x,y,\xi)u(y) dy d\xi,
\\
&a(x,y,\xi)=a_1(x,\xi)a_2(y)\quad\text{or}\quad a(x,y,\xi)=a_2(x)a_1(y,\xi),
\end{align*}
where
$a_1\in C^\infty\p{\R^n\times\R^n}$, $a_2\in L^\infty\p{\R^n}$,
and $\varphi(y,\xi)\in C^\infty\p{\R^n_y\times\R^n_\xi}$ is a real-valued
function.
Assume that
\[
\abs{\partial_x^\alpha\partial_\xi^\beta a_1(x,\xi)}
\leq C_{\alpha\beta}
\]
for $|\alpha|,|\beta|\leq 2n+1$.
Also assume that, on $\R^n\times \tilde{E}$,
\[
\abs{\det\partial_y\partial_{\xi}\varphi(y,\xi)}\geq C>0,
\]
and each entry $h(y,\xi)$ of
the matrix $\partial_y\partial_\xi\varphi(y,\xi)$ satisfies
\[
\abs{\partial_y^\alpha h(y,\xi)}\leq C_{\alpha},
\qquad
\abs{\partial_\xi^\beta h(y,\xi)}\leq C_{\beta}
\]
for $|\alpha|,|\beta|\leq 2n+1$, where
\[
\tilde{E}=\bigcup_{x,y\in\R^n}E_{x,y}\,;\quad E_{x,y}=\supp a(x,y,\cdot).
\]
Then $T$ is $L^2(\R^n)$-bounded,
and satisfies
\[
\n{T}_{L^2\to L^2}\leq C
\sup_{|\alpha|,|\beta|\leq 2n+1}
 \n{\partial_y^\alpha\partial_\xi^\beta a(y,\xi)}
  _{L^\infty\p{\R^n_y\times\R^n_\xi}}.
\]
\end{cor}
\medskip
\begin{proof}
Note that $T$ is a product of the multiplication of the function $a_2$
and the operator defined by \eqref{T1} or \eqref{T2}, which are all
$L^2$-bounded by the assumption.
\end{proof}
\medskip
\par
%\newpage
\bigskip
\section{Weighted $L^2$-estimates}
Asada and Fujiwara \cite{AF} proved Corollary \ref{Cor2.6} without the
product type
assumption for $a(x,y,\xi)$, but assumed the boundedness of all the derivatives
of $a(x,y,\xi)$ and that of each entry of the matrix
$\partial_\xi\partial_\xi\varphi$.
The following theorem,
which is a generalized version of Theorem \ref{mainth},
says that we do not need the boundedness assumption for
$\partial_\xi\partial_\xi\varphi$ if $a(x,y,\xi)$ has a decaying property.
In this case, we have weighted estimates as follows.
For $m\in\R$, we use the notation
\[
\langle x\rangle^m=\p{1+|x|^2}^{m/2},
\]
and let $L^2_m(\R^n)$ be the set of functions $f$ such that the norm
\[
\n{f}_{L^2_m(\R^n)}
=\p{\int_{\R^n}\abs{\langle x\rangle^m f(x)}^2\,dx}^{1/2}
\]
is finite.
\medskip
\begin{theorem}\label{Th3.1}
Suppose $m_1,m_2\in\R$.
Let the operator $T$ be defined by
\begin{equation}\label{T3}
 Tu(x)
 =\int_{\R^n}\int_{\R^n} e^{i(x\cdot\xi+\varphi(y,\xi))}a(x,y,\xi)u(y) dy d\xi,
\end{equation}
where $a(x,y,\xi)\in C^\infty\p{\R^n_x\times\R^n_y\times\R^n_\xi}$, and
$\varphi(y,\xi)\in C^\infty\p{\R^n_y\times\R^n_\xi}$ is a real-valued function.
Assume that
\[
\abs{\det\partial_y\partial_{\xi}\varphi(y,\xi)}\geq C>0\quad
\text{on $\R^n\times \tilde{E}$},
\]
where
\[
\tilde{E}=\bigcup_{x,y\in\R^n}E_{x,y}\,;\quad E_{x,y}=\supp a(x,y,\cdot).
\]
Also assume one of the followings:
\begin{enumerate}
\item
For all $\alpha$, $\beta$, and $\gamma$,
\[
\abs{\partial_x^\alpha\partial_y^\beta\partial_\xi^\gamma
a(x,y,\xi)} \leq C_{\alpha\beta\gamma}\langle
x\rangle^{m_1-|\alpha|}\jp{y}^{m_2},
\]
and for all $|\alpha|\geq1$ and $|\beta|\geq1$,
\[
\abs{\partial_\xi^\beta\varphi(y,\xi)}\leq C_{\beta}\langle y\rangle,
\quad
\abs{\partial_y^\alpha\partial_\xi^\beta\varphi(y,\xi)} \leq
C_{\alpha\beta}
\quad
\text{on}\quad \R^n\times \tilde{E}.
\]
\item
For all $\alpha$, $\beta$, and $\gamma$,
\[
\abs{\partial_x^\alpha\partial_y^\beta\partial_\xi^\gamma
a(x,y,\xi)} \leq C_{\alpha\beta\gamma}\jp{x}^{m_1}\langle
y\rangle^{m_2-|\beta|},
\]
and for all $\alpha$ and $|\beta|\geq1$,
\[
\abs{\partial_y^\alpha\partial_\xi^\beta\varphi(y,\xi)} \leq
C_{\alpha\beta}\langle y\rangle^{1-|\alpha|}
\quad
\text{on}\quad \R^n\times \tilde{E}.
\]
\end{enumerate}
Then $T$ is bounded from $L^2_{m+m_1+m_2}(\R^n)$ to $L^2_{m}(\R^n)$
for any $m\in\R$.
\end{theorem}
\medskip
\begin{rem}
From the assumptions for phase functions $\varphi$ in Theorem \ref{Th3.1},
we obtain the estimate
\begin{equation}\label{below}
C_1\langle y\rangle\leq \langle \partial_\xi\varphi(y,\xi)\rangle
\leq C_2\langle y\rangle
\quad\text{on $\R^n\times \tilde{E}$}
\end{equation}
for some $C_1,C_2>0$.
In fact, the estimate $\jp{\partial_\xi\varphi(y,\xi)}\leq C_2\jp{y}$
is obtained from any assumptions (1) or (2).
Especially, we have $\jp{\partial_\xi\varphi(0,\xi)}\leq C$.
From the expression
\[
\partial_\xi\varphi(y,\xi)-\partial_\xi\varphi(0,\xi)
=\partial_y\partial_\xi\varphi(z,\xi)y
\]
with some $z\in\R^n$,
we obtain
\begin{align*}
|y|&\leq C\abs{\partial_\xi\varphi(y,\xi)-\partial_\xi\varphi(0,\xi)}
\\
&\leq C\abs{\partial_\xi\varphi(y,\xi)}+C\abs{\partial_\xi\varphi(0,\xi)}
\end{align*}
by the assumptions for $\varphi$.
Hence we have the estimate
$C_1\langle y\rangle\leq \langle \partial_\xi\varphi(y,\xi)\rangle$,
as well.
\end{rem}
\medskip
\begin{proof}
We show the $L^2$-boundedness of the operator $T_b$ defined by
\[
 T_bu(x)
 =\int\int e^{i(x\cdot\xi+\varphi(y,\xi))}b(x,y,\xi)u(y) dy d\xi,
\]
where
\[
 b(x,y,\xi)=\langle x\rangle^{m} a(x,y,\xi) \langle y\rangle^{-(m+m_1+m_2)}.
\]
By using the cut-off function $\chi(x)\in C^\infty_0\p{|x|\leq 1/2}$ which is
equal to one near the origin, we decompose $b$ into two parts:
\begin{align*}
&b^I(x,y,\xi)
=b(x,y,\xi)
 \chi\p{\p{x+\partial_\xi\varphi(y,\xi)}
     /\langle \partial_\xi\varphi(y,\xi)\rangle},
\\
&b^{II}(x,y,\xi)
=b(x,y,\xi)
  \p{1-\chi}\p{\p{x+\partial_\xi\varphi(y,\xi)}
     /\langle \partial_\xi\varphi(y,\xi)\rangle}.
\end{align*}
The corresponding decomposition of the operator $T_b$ is denoted by
$T^I$ and $T^{II}$ respectively.
\par
On the support of $b^I(x,y,\xi)$, we have
$\abs{x+\partial_\xi\varphi(y,\xi)}
\leq(1/2) \langle \partial_\xi\varphi(y,\xi)\rangle$,
hence we have the estimates
\[
|x|\leq
  \abs{\partial_\xi\varphi(y,\xi)}
 +\frac12 \langle \partial_\xi\varphi(y,\xi)\rangle,
\qquad
\abs{\partial_\xi\varphi(y,\xi)}
\leq |x|+\frac12 \langle \partial_\xi\varphi(y,\xi)\rangle.
\]
From the first estimate and estimate \eqref{below}, we obtain
$\langle x\rangle\leq C\langle y\rangle$.
From the second estimate, we obtain
$\langle \partial_\xi\varphi(y,\xi)\rangle
\leq 2\langle x\rangle+(1/2) \langle \partial_\xi\varphi(y,\xi)\rangle$,
hence
$\langle \partial_\xi\varphi(y,\xi)\rangle\leq 4\langle x\rangle$,
which implies
$\langle y\rangle\leq C\langle x\rangle$
by \eqref{below} again.
Thus we have the equivalence of $\langle y\rangle$ and $\langle x\rangle$,
and obtain
\begin{equation}\label{xdecay}
\abs{\partial_x^\alpha\partial_y^\beta\partial_\xi^\gamma b^I(x,y,\xi)}
\leq C_{\alpha\beta\gamma}\langle x\rangle^{-|\alpha|}
\end{equation}
or
\begin{equation}\label{ydecay}
\abs{\partial_x^\alpha\partial_y^\beta\partial_\xi^\gamma b^I(x,y,\xi)}
\leq C_{\alpha\beta\gamma}\langle y\rangle^{-|\beta|}
\end{equation}
from the assumptions (1) and (2) respectively.
\par
We assume estimate \eqref{xdecay}. Otherwise, assume
\eqref{ydecay} and just change the role of $x$ and $y$ below. Let
real-valued positive functions $\Phi_0(x)$,
$\Phi_k(x)=\Phi(x/2^k)$ ($k\in\N$) form a partition of unity which
satisfy $\supp \Phi_0\subset \{x;|x|< 2\}$, $\supp \Phi\subset
\{x;1/2<|x|<2\}$. We decompose $b^I$ into the sum of
$b^I_k(x,y,\xi)=\Phi_k(x)b^I(x,y,\xi)$. By the equivalence of
$\langle x\rangle$ and $\langle y\rangle$ on the support of $b^I$,
we can write
\[
b^I_k(x,y,\xi)
=\Phi_k(x)b^I(x,y,\xi)\tilde{\Psi}_k(y)
\]
with functions $\tilde{\Psi}_k\in C_0^\infty(\R^n)$ which are of the
from $\tilde{\Psi}_k(y)=\tilde{\Psi}(y/2^k)$,
$\tilde{\Psi}\in C_0^\infty(\R^n\setminus0)$
with large $k$.
Furthermore, we have
\[
b^I_k(2^kx,y,\xi)
=\Psi_k(2^kx)\sum_{l\in\Z^n}e^{il\cdot x}b_{kl}(y,\xi)\tilde{\Psi}_k(y),
\]
where $\Psi_k$ is the characteristic function of the support of $\Phi_k$, and
\begin{align*}
b_{kl}(y,\xi)
&=\int e^{-il\cdot x} b^I_k(2^kx,y,\xi)\,dx
\\
&=\p{1+|l|^2}^{-n}\int e^{-il\cdot x} \p{1-\triangle_x}^n
    \b{\Phi_k(2^kx)b^I(2^kx,y,\xi)}\,dx
\end{align*}
is the Fourier coefficients of the function
$b^I_k(2^kx,y,\xi)$ in the variable $x$.
Then, by estimate \eqref{xdecay}, we have
\[
\abs{\partial_y^\alpha\partial_\xi^\beta b_{kl}(y,\xi)}
\leq C_{\alpha\beta}\p{1+|l|^2}^{-n},
\]
where $C_{\alpha\beta}$ is independent of $k,l\in\Z^n$.
Thus we have the decomposition
\[
T^I=\sum_{l\in \Z^n}
  \sum_{k\in \Z^n}e^{il\cdot x/2^k}\Psi_kT_{kl}\tilde{\Psi}_k,
\]
where
\[
T_{kl}v(x)=
\int\int
 e^{i(x\cdot\xi+\varphi(y,\xi))}b_{kl}(y,\xi)v(y) dyd\xi.
\]
We remark that
\begin{align*}
\n{\sum_{k\in \Z^n}e^{il\cdot x/2^k}\Psi_kT_{kl}\tilde{\Psi}_ku}_{L^2}^2
&\leq C\sum_{k\in \Z^n}\n{\Psi_kT_{kl}\tilde{\Psi}_ku}_{L^2}^2
\\
& \leq C\sup_{k\in \Z^n}\n{T_{kl}}_{L^2\to L^2}^2
  \sum_{k\in \Z^n}\n{\tilde{\Psi}_ku}_{L^2}^2
\\
& \leq C\p{1+|l|^2}^{-2n}\n{u}_{L^2}^2
\end{align*}
by Corollary \ref{Cor2.6}.
Hence we have
\begin{align*}
\n{T^I}_{L^2\to L^2}
&\leq C\sum_{l\in \Z^n}\p{1+|l|^2}^{-n}
\\
&\leq C,
\end{align*}
that is, the $L^2$-boundedness of $T^I$.
\par
Next, we show the boundedness of $T^{II}$.
Let $\rho\in C^\infty_0$ be a real-valued function which satisfies
\[
\sum_{k\in \Z^n}\rho(\xi-k)=1.
\]
We decompose $b^{II}(x,y,\xi)$ into the sum of
\[
b^{II}_k(x,y,\xi)=b^{II}(x,y,\xi)\rho(\xi-k)
\]
and set
\[
 T_ku(x)
 =\int_{\R^n}\int_{\R^n} e^{i(x\cdot\xi+\varphi(y,\xi))}
  b^{II}_k(x,y,\xi)u(y) dy d\xi.
\]
We claim, we may replace $b^{II}_k(x,y,\xi)$ by the symbol (denoted
by $b^{II}_k(x,y,\xi)$ again) which has the same (or smaller) support
and satisfies
the estimate
\begin{equation}\label{b_k^II}
\abs{\partial_x^\alpha\partial_y^\beta\partial_\xi^\gamma b^{II}_k(x,y,\xi)}
\leq C_{\alpha\beta\gamma}\langle x\rangle^{-(n+1)}\langle y\rangle^{-(n+1)},
\end{equation}
where $C_{\alpha\beta\gamma}$ is independent of $k\in\Z^n$.
Indeed, by integration by parts we have
\[
 T_ku(x)
 =\int\int e^{i(x\cdot\xi+\varphi(y,\xi))}
  L^Nb^{II}_k(x,y,\xi)u(y) dy d\xi,
\]
where $L$ is the transpose of the operator
\[
^tL=\frac{x+\partial_\xi\varphi}{i\abs{x+\partial_\xi\varphi}^2}
    \cdot \partial_\xi
\]
and $N$ is a positive integer.
We have
$\langle \partial_\xi\varphi(y,\xi)\rangle
\leq C\abs{x+\partial_\xi\varphi(y,\xi)}$
on the support of $b^{II}(x,y,\xi)$,
hence we have
\begin{align*}
 &\langle x\rangle
\leq
  \abs{x+\partial_\xi\varphi(y,\xi)}
 +2\langle \partial_\xi\varphi(y,\xi)\rangle
\leq C\abs{x+\partial_\xi\varphi(y,\xi)},
\\
 &\langle y\rangle
\leq C\langle \partial_\xi\varphi(y,\xi)\rangle
\leq C\abs{x+\partial_\xi\varphi(y,\xi)}
\end{align*}
by estimate \eqref{below}.
Thus
$\abs{x+\partial_\xi\varphi}^{-1}$ is dominated by $\langle
x\rangle^{-1}$ and $\langle y\rangle^{-1}$, and we can justify our
claim by taking large $N$.
\par
Let $T_k^*$ be the adjoint of $T_k$, and we have
\[
T_kT_l^*v(x)=\int K_{kl}(x,y)v(y)\,dy,\qquad
T_k^*T_lv(x)=\int \tilde{K}_{kl}(x,y)v(y)\,dy,
\]
where
\begin{align*}
&K_{kl}(x,y)
=\int\int\int e^{i\b{x\cdot\xi-y\cdot\eta+\varphi(z,\xi)-\varphi(z,\eta)}}
 b^{II}_k(x,z,\xi)\overline{b^{II}_l(y,z,\eta)}
 \,dz d\xi d\eta,
\\
&\tilde{K}_{kl}(x,y)
=\int\int\int e^{i\b{\varphi(y,\xi)-\varphi(x,\eta)+z\cdot(\xi-\eta)}}
 b^{II}_l(z,y,\xi)\overline{b^{II}_k(z,x,\eta)}
 \,dz d\xi d\eta.
\end{align*}
By integration by parts, we have
\begin{align*}
&\int e^{i\p{\varphi(z,\xi)-\varphi(z,\eta)}}
b^{II}_k(x,z,\xi)\overline{b^{II}_l(y,z,\eta)}\,dz
\\
=&\int e^{i\p{\varphi(z,\xi)-\varphi(z,\eta)}}
L^{2n+1}\p{b^{II}_k(x,z,\xi)\overline{b^{II}_l(y,z,\eta)}}\,dz,
\end{align*}
where $L$ is the transpose of the operator
\[
^tL=\frac 1i
\frac
{\partial_z\varphi(z,\xi)-\partial_z\varphi(z,\eta)}
{\abs{\partial_z\varphi(z,\xi)-\partial_z\varphi(z,\eta)}^2}
\cdot\partial_z.
\]
From the assumptions, we obtain
\[
\abs{\partial_z\varphi(z,\xi)-\partial_z\varphi(z,\eta)}\geq C|\xi-\eta|
\]
and
\[
\abs{\partial_z^\beta\varphi(z,\xi)-\partial_z^\beta\varphi(z,\eta)}
\leq C_\beta|\xi-\eta|
\]
for all $\beta$.
From this argument and \eqref{b_k^II}, we obtain
\[
|K_{kl}(x,y)|
\leq C \langle x\rangle^{-(n+1)}
       \langle y\rangle^{-(n+1)}
  \p{1+|k-l|^{2n+1}}^{-1},
\]
where $C$ is independent of $k,l\in\Z^n$.
Then we have
\begin{align*}
&\sup_{x}\int |K_{kl}(x,y)|\,dy
\leq C \p{1+|k-l|^{2n+1}}^{-1},
\\
&\sup_{y}\int |K_{kl}(x,y)|\,dx
\leq C \p{1+|k-l|^{2n+1}}^{-1}
\end{align*}
which implies, by Lemma \ref{Lem2.1},
\[
\n{T_kT_l^*}_{L^2\to L^2}
\leq C \p{1+|k-l|^{2n+1}}^{-1}.
\]
Similarly, we have
\[
\n{T_k^*T_l}_{L^2\to L^2}
\leq C \p{1+|k-l|^{2n+1}}^{-1}
\]
if we take
\[
^tL=\frac 1i
\frac
{\xi-\eta}{\abs{\xi-\eta}^2}
\cdot\partial_z.
\]
Then we have
\[
\n{T_kT_l^*}_{L^2\to L^2},\,
\n{T_k^*T_l}_{L^2\to L^2}
    \leq C \b{\gamma\p{k-l}}^2,
\]
where
\[
\gamma\p{j}
=\p{1+|j|^{2n+1}}^{-1/2}
\]
and it satisfies the estimate
\[
\sum_{j\in\Z^n}
\gamma\p{j}<\infty.
\]
By Lemma \ref{Lem2.2}, we have the $L^2$-boundedness of $T^{II}$.
\end{proof}
\medskip
\par
%\newpage
\bigskip
\section{Applications}
In this section, we explain how to use Theorem \ref{Th3.1} to show the
smoothing effect of generalized Schr\"odinger equations.
The main tool is a class of Fourier integral operators
of the form
\begin{equation}\label{eq:defI}
\begin{aligned}
&T_\psi u(x)
=(2\pi)^{-n}\int_\Rn\int_\Rn e^{i(x\cdot\xi-y\cdot\psi(\xi))}u(y) dyd\xi,
\\
&T_\psi^{-1}u(x)
=(2\pi)^{-n}\int_\Rn\int_\Rn e^{i(x\cdot\xi-y\cdot\psi^{-1}(\xi))}
 u(y) dyd\xi,
\end{aligned}
\end{equation}
where $\psi,\psi^{-1}:\Rn\setminus0\to\Rn\setminus0$ are $C^\infty$-maps
satisfying $\psi\circ\psi^{-1}(\xi)=\psi^{-1}\circ\psi(\xi)=\xi$,
$\psi(\lambda\xi)=\lambda\psi(\xi)$, and
$\psi^{-1}(\lambda\xi)=\lambda\psi^{-1}(\xi)$
for all $\lambda>0$ and $\xi\in\Rn\setminus0$.
We remark that we have
\begin{equation}\label{Fourier}
T_\psi u(x)=F^{-1}_\xi\br{\p{F_xu}(\psi(\xi))}(x),\qquad
T_\psi^{-1}u(x)=F^{-1}_\xi\br{\p{F_xu}(\psi^{-1}(\xi))}(x),
\end{equation}
where $F_x$ ($F_\xi^{-1}$ {\it resp.})
denotes the (inverse {\it resp.})
Fourier transform.
Hence, we have $T_\psi^{-1}\cdot T_\psi =T_\psi\cdot T_\psi^{-1}=id$,
and the formula
\begin{equation}\label{eq:transg}
T_\psi\cdot a(D)\cdot T_\psi^{-1}=\p{a\circ \psi}(D),
\end{equation}
where $a(D)=F^{-1}_\xi a(\xi)F_x$.
By \eqref{Fourier} and Plancherel's theorem,
the operators $T_\psi$ and $T_\psi^{-1}$ are $L^2$-bounded.
Furthermore, as a corollary of Theorem \ref{Th3.1}, we have the following:
\medskip
\begin{cor}\label{Cor4.1}
Suppose $m\in\Z$ and $|m|<n/2$.
Assume that $\abs{\det \partial\psi(\xi)}\geq C>0$.
Then the operators $T_\psi$ and $T_\psi^{-1}$ defined by
\eqref{eq:defI} are $L^2_{m}(\R^n)$-bounded.
\end{cor}
\medskip
\begin{rem}\label{Rem4.1}
By Corollary \ref{Cor4.1} and the interpolation,
we have the $L^2_{m}(\R^n)$-boundedness of $T_\psi$
and $T_\psi^{-1}$ with $m\in\R$ such that $|m|\leq[n/2]_{-}$
($[k]_{-}$ denotes the greatest integer less than $k$).
\end{rem}
\medskip
\begin{proof}
We prove the boundedness of $T_\psi$, from which the boundedness of
$T_\psi^{-1}$ follows.
Let $\chi(\xi)\in C_0^\infty$ be a cut off function of the origin.
By \eqref{Fourier}, we have
\[
\p{1-\chi(D)}T_\psi u(x)
=(2\pi)^{-n}\int_\Rn\int_\Rn e^{i(x\cdot\xi-y\cdot\psi(\xi))}
 \p{1-\chi(\xi)}u(y) dyd\xi.
\]
Since $\psi(\xi)$ is smooth away from the origin,
$\p{1-\chi(D)}T_\psi$ is $L^2_m$-bounded by Theorem \ref{Th3.1}.
On the other hand, if we note
\[
e^{ix\cdot\xi}=\frac{1-ix\cdot\partial_\xi}{\jp{x}^2}e^{ix\cdot\xi},
\qquad
e^{-iy\cdot\xi}=\frac{1+iy\cdot\partial_\xi}{\jp{y}^2}e^{-iy\cdot\xi},
\]
we have, by change of variables and integration by parts,
\begin{equation}\label{plus}
\begin{aligned}
&\chi(D)T_\psi u(x)
\\
&=(2\pi)^{-n}\int\int e^{i(x\cdot\xi-y\cdot\psi(\xi))}\chi(\xi) u(y) dyd\xi
\\
&=(2\pi)^{-n}\int\int e^{i(x\cdot\xi-y\cdot\psi(\xi))}
\p{
\frac{1+ix\cdot\partial\chi(\xi)+x\chi(\xi){}^t\partial\psi(\xi){}^ty}
     {\jp{x}^2}
}
u(y) dyd\xi
\\
&=\frac1{\jp{x}^2}T_\psi u
+i\frac{x}{\jp{x}^2}\cdot\partial\chi(D)T_\psi u
+\frac{x}{\jp{x}^2}\chi(D){}^t\partial\psi(D)T_\psi \p{{}^txu}
\end{aligned}
\end{equation}
and
\begin{equation}\label{minus}
\begin{aligned}
\chi(D)T_\psi u(x)
&=(2\pi)^{-n}\int\int e^{i(x\cdot\psi^{-1}(\xi)-y\cdot\xi)}
 \chi\p{\psi^{-1}(\xi)}\abs{\det\partial\psi^{-1}(\xi)} u(y) dyd\xi
\\
&=(2\pi)^{-n}\int\int e^{i(x\cdot\psi^{-1}(\xi)-y\cdot\xi)}
\p{
\frac{1+{\bf a}(\xi)\cdot y+xA(\xi){}^t y}{\jp{y}^2}
}
u(y) dyd\xi
\\
&=d(D)T_\psi\p{\frac{u}{\jp{x}^2}}
+|D|^{-1}|D|{\bf a}\p{\psi(D)}d(D)\cdot T_\psi\p{\frac{x}{\jp{x}^2}u}
\\
&+xA(\psi(D))d(D) T_\psi\p{\frac{{}^t x}{\jp{x}^2}u}
\end{aligned}
\end{equation}
where
\begin{align*}
&A(\xi)=\chi\p{\psi^{-1}(\xi)}
        \abs{\det\partial\psi^{-1}(\xi)}\partial\psi^{-1}(\xi),
\quad
{\bf a}(\xi)
 =-i\partial\b{\chi\p{\psi^{-1}(\xi)}\abs{\det\partial\psi^{-1}(\xi)}},
\\
&d(\xi)=\abs{\det\partial\psi(\xi)}.
\end{align*}
We remember here that $T_\psi$ is $L^2$-bounded.
Assume that $T_\psi$ is $L^2_{\pm(k-1)}$-bounded with some $k<n/2$, $k\in\N$.
We remark that $\chi(D)$, $d(D)$ and all entries of $\partial\chi(D)$,
$\partial\psi(D)$, $A(\psi(D))$, $|D|{\bf a}\p{\psi(D)}$
are $L^2_{\pm (k-1)}$-bounded, and $|D|^{-1}$ is bounded from
$L^2_{-(k-1)}$ to $L^2_{-k}$.
To justify these boundedness, use the results of Kurtz and Wheeden \cite{KW},
Stein and Weiss \cite{SW}.
Using them, we obtain the $L^2_{k}$-boundedness of $\chi(D)T_\psi$ from
\eqref{plus}, and $L^2_{-k}$-boundedness from \eqref{minus}.
Then, by induction, we have the desired result.
\end{proof}
\medskip
Now, let $p(\xi)\in C^\infty(\R^n\setminus0)$ be a positive function which
satisfies $p(\lambda\xi)=\lambda p(\xi)$
for all $\lambda>0$ and $\xi\in\Rn\setminus0$,
and let
\[
L_p=p(D_x)^2=F^{-1}_\xi p(\xi)^2F_x
\]
be the corresponding Fourier multiplier.
Assume that $\Sigma=\b{\xi;p(\xi)=1}$ has non-vanishing Gaussian curvature.
We consider a generalized Schr\"odinger equation
\begin{equation}
\label{eq:GSch}
\left\{
\begin{array}{l}
\p{i\partial_t+L_p}u(t,x)=0, \\
u(0,x)=f(x).
\end{array}
\right.
\end{equation}
If we take
\begin{equation}\label{eq:phase}
\psi(\xi)=p(\xi)\frac{\nabla p(\xi)}{\abs{\nabla p(\xi)}},
\end{equation}
we have the relation
\begin{equation}\label{eq:trans}
T_\psi\cdot (-\triangle_x)\cdot T_\psi^{-1}=L_p
\end{equation}
by \eqref{eq:transg}, and the $L^2_{-1}$-boundedness of
the operators $T_\psi$ and $T_\psi^{-1}$ by Corollary \ref{Cor4.1}.
In fact, the curvature condition on $\Sigma$ means that
the Gauss map
\[
\frac{\nabla p}{\abs{\nabla p}}:\Sigma\to S^{n-1}
\]
is a global diffeomorphism and its Jacobian never vanishes
(see Kobayashi and Nomizu \cite{KN}).
Hence, we can construct the inverse $C^\infty$-map $\psi^{-1}(\xi)$
of $\psi(\xi)$, and can justify the assumption of Corollary \ref{Cor4.1}.
Applying $T_\psi^{-1}$ defined by \eqref{eq:defI}
with \eqref{eq:phase} to equation
(\ref{eq:GSch}),
and introducing $v=T_\psi^{-1}u$ and $g=T_\psi^{-1}f$,
\eqref{eq:GSch} can be transformed to the equation
\begin{equation}
\label{eq:OSchc}
\left\{
\begin{array}{l}
\p{i\partial_t-\Delta_x} v(t,x)=0,  \\
v(0,x)=g(x),
\end{array}
\right.
\end{equation}
by \eqref{eq:trans}.
It has been already known that
classical Schr\"odinger equation \eqref{eq:OSchc} has the global smoothing
estimate
\begin{equation}\label{eq:sm}
\n{\sigma(X,D) v}_{L^2\p{\R_t\times\R^n_x}}\leq C\n{g}_{L^2\p{\R^n_x}},
\end{equation}
where $n\geq3$ and 
\[
\sigma(X,D)=\langle x\rangle^{-1} \langle D\rangle^{1/2}.
\]
See Ben-Artzi and Klainerman \cite{BK}, Simon \cite{Si},
 Kato and Yajima \cite{KY}, or Walther \cite{Wa}.
From this fact, we can extract a similar estimate
for generalized Schr\"odinger equation (\ref{eq:GSch}).
In fact, we have
\[
\langle D\rangle^{1/2}u
= M\p{1+p(D)^2}^{1/4}T_\varphi v
= M T_\varphi \langle D\rangle^{1/2} v
\]
where
\[
M= \langle D\rangle^{1/2}\p{1+p(D)^2}^{-1/4}.
\]
Here we have used the formula \eqref{eq:transg} with
$a(\xi)=\p{1+|\xi|^2}^{1/4}$.
Hence we have
\[
\sigma(X,D)u = \langle x\rangle^{-1}M T_\psi\langle
x\rangle\sigma(X,D)v.
\]
Since $M$ is $L^2_{-1}(\Rnx)$-bounded
by Theorem \ref{mainth}, and $T_\psi$ by Corollary \ref{Cor4.1},
we obtain
\[
\n{\sigma(X,D) u}_{L^2\p{\R_t\times\R^n_x}} \leq
C\n{f}_{L^2\p{\R^n_x}}
\]
from estimates \eqref{eq:sm} and
\[
\n{g}_{L^2\p{\R^n_x}}=\n{T_\psi^{-1}f}_{L^2\p{\R^n_x}}
\leq C\n{f}_{L^2\p{\R^n_x}}.
\]
\par
Thus, we have obtained the following result which was partially
proved for a type of
polynomial $p(\xi)^2$ by Ben-Artzi and Devinatz \cite{BD},
and fully for radially symmetric $p(\xi)^2$ by Walther \cite{Wa2}.
\medskip
\begin{theorem}\label{Th4.2}
Suppose $n\geq3$.
Assume that $\Sigma=\b{\xi;p(\xi)=1}$ has non-vanishing Gaussian curvature.
Then the solution $u(t,x)$ to equation \eqref{eq:GSch} has the estimate
\[
\n{\langle x\rangle^{-1}
 \langle D\rangle^{1/2} u}_{L^2\p{\R_t\times\R^n_x}}
\leq C\n{f}_{L^2\p{\R^n_x}}.
\]
\end{theorem}
\medskip
In Theorem \ref{Th4.2}, the order $\lq\lq-1"$
for the weight is the best possible
one because of the estimate for the low frequency part
(Walther \cite{Wa}, \cite{Wa2}).
But, if we replace $\jp{D}^{1/2}$ by $|D|^{1/2}$,
we have another type of estimate
\[
\n{\langle x\rangle^{-\delta}
 |D|^{1/2} u}_{L^2\p{\R_t\times\R^n_x}}
\leq C\n{f}_{L^2\p{\R^n_x}}
\]
for $\delta>1/2$.
Chihara \cite{Chih} obtained this type of estimates for rather general
$p(\xi)^2$.
In our forthcoming paper \cite{RS3}, we use our main result
Theorem \ref{Th3.1} to obtain a refinement of this estimate.
\par
%%%%%%%%%%%%%%%%%%%%%%%%%%%%%references%%%%%%%%%%%%%%%%%%%%%%%%%%%%%%%%%%%
%\newpage

\end{document}